\documentclass[12pt,a4paper]{article}

\newtheorem{theorem}{Theorem}[section]
\newtheorem{proposition}[theorem]{Proposition}

\newtheorem{corollary}[theorem]{Corollary}
\newtheorem{example}[theorem]{Example}

\usepackage{hyperref}
\usepackage{amsmath}
\usepackage{amssymb}
\usepackage{array}
\usepackage[english]{babel}
\usepackage{calc}
\usepackage{enumerate}
\usepackage[latin1]{inputenc}
\usepackage{latexsym}

\newcommand{\Rp}{\mathbb R^+}

\newcommand{\Lom}{L^1(\omega)}
\newcommand{\Loo}{L^1(\omega_1)}
\newcommand{\Lot}{L^1(\omega_2)}
\newcommand{\Mom}{M(\omega)}
\newcommand{\Moo}{M(\omega_1)}
\newcommand{\Mot}{M(\omega_2)}
\newcommand{\Com}{C_0(1/\omega)}
\newcommand{\Coo}{C_0(1/\omega_1)}
\newcommand{\Cot}{C_0(1/\omega_2)}
\newcommand{\Lim}{L^{\infty}(1/\omega)}
\newcommand{\Lio}{L^{\infty}(1/\omega_1)}
\newcommand{\Lit}{L^{\infty}(1/\omega_2)}

\begin{document}
\title{Generalized Hardy-Ces\`aro operators between weighted spaces}
\author{Thomas Vils Pedersen}

%\date
\maketitle

\footnotetext{2010 {\em Mathematics Subject Classification:} 
44A15, 47B34, 47B38, 47G10}
\footnotetext{{\em Keywords:} 
Generalized Hardy-Ces\`aro operators, weighted spaces, weak compactness.}

\begin{abstract}
\noindent
We characterize those non-negative, measurable functions $\psi$ on $[0,1]$ and 
positive, continuous functions $\omega_1$ and $\omega_2$ on $\Rp$ for which
the generalized Hardy-Ces\`aro operator
$$(U_{\psi}f)(x)=\int_0^1 f(tx)\psi(t)\,dt$$
defines a bounded operator $U_{\psi}:\Loo\to\Lot$. Furthermore, we extend
$U_{\psi}$ to a bounded operator on $\Moo$ with range in 
$\Lot\oplus\mathbb C\delta_0$.
Finally, we show that the zero operator is the only weakly compact
generalized Hardy-Ces\`aro operator from $\Loo$ to $\Lot$.
\end{abstract}

\section{Introduction}
\label{sec:intro}

A classical result of Hardy (\cite{Ha-Li-Po-2}) shows that the Hardy-Ces\`aro 
operator
$$(Uf)(x)=\frac{1}{x}\int_0^x f(s)\,ds$$
defines a bounded linear operator on $L^p(\mathbb R^+)$ with
$\|U\|=p/(p-1)$ for $p>1$. Clearly, $U$ is not bounded on $L^1(\mathbb R^+)$. 
Hardy's result has been generalized in various ways, of which we will mention
some, which have inspired this paper.

For $1\le p\le q\le\infty$ and 
non-negative measurable functions $u$ and $v$ on $\Rp$, 
Muckenhoupt (\cite{Muc}) and Bradley (\cite{Bra}) gave a necessary
and sufficient condition for the existence of a constant $C$ such that
$$\left(\int_0^{\infty}\left(u(x)\int_0^xf(t)\,dt\right)^q\,dx\right)^{1/q}
\le C\left(\int_0^{\infty}\left(v(x)f(x)\right)^p\,dx\right)^{1/p}$$
for every positive, measurable function $f$ on $\Rp$. This can be rephrased 
as a characterization of the weighted $L^p$ and $L^q$ spaces on $\Rp$ between which
the Hardy-Ces\`aro operator $U$ is bounded.

In a different direction, for a non-negative measurable funtion $\psi$ on 
$[0,1]$, Xiao (\cite{Xi}) considered the 
generalized Hardy-Ces\`aro operators
$$(U_{\psi}f)(x)=\int_0^1 f(tx)\psi(t)\,dt=\frac{1}{x}\int_0^x f(s)\psi(s/x)\,ds$$
for measurable functions $f$ on $\mathbb R^n$. Xiao proved that 
$U_{\psi}$ defines a bounded operator on $L^p(\mathbb R^n)$ (for $p\ge1$) 
if and only if
$$\int_0^1\frac{\psi(t)}{t^{n/p}}\,dt<\infty.$$

Finally, we mention that Albanese, Bonet and Ricker in a recent series
of papers (see, for instance, \cite{Al-Bo-Ri}) have considered the spectrum, 
compactness and other properties of the Hardy-Ces\`aro operator on
various function spaces.

In this paper we will study the generalized Hardy-Ces\`aro operators between weighted
spaces of integrable functions. Let $\omega$ be a positive, continuous function on $\Rp$
and let $\Lom$ be the Banach space of (equivalence classes of) measurable functions $f$ on 
$\mathbb R^+$ for which 
$$\|f\|_{\Lom}=\int_0^{\infty}|f(t)|\omega(t)\,dt<\infty.$$
In the usual way we identify the dual space of $\Lom$ with the space $\Lim$ of
measurable functions $h$ on $\Rp$ for which
$$\|h\|_{\Lim}=\text{ess\,sup}_{t\in\Rp}|h(t)|/\omega(t)<\infty.$$
We denote by $\Com$ the closed subspace of $\Lim$ consisting of the continuous functions $g$ in 
$\Lim$ for which $g/\omega$ vanishes at infinity. 
Finally, we identify the dual space of $\Com$ with the space $\Mom$ of
locally finite, complex Borel measures $\mu$ on $\mathbb R^+$ for which 
$$\|\mu\|_{\Mom}=\int_{\Rp}\omega(t)\,d|\mu|(t)<\infty.$$
We consider the space $\Lom$ as a closed subspace of $\Mom$.

In Section~\ref{sec:char} we characterize those functions $\psi, \omega_1$ and 
$\omega_2$ for which $U_{\psi}$ defines a bounded operator from $\Loo$ to 
$\Lot$. These operators are extended to bounded operators on $\Moo$ in 
Section~\ref{sec:ext}, where we also obtain results about their ranges. 
Finally, in Section~\ref{sec:wk} we show that there are no non-zero 
weakly compact generalized Hardy-Ces\`aro operators from $\Loo$ to 
$\Lot$.

\section{A characterization of the generalized Hardy-Ces\`aro operators}
\label{sec:char}

For a non-negative, measurable function $\psi$ on $[0,1]$ and positive, 
continuous functions $\omega_1$ and $\omega_2$ on $\Rp$, 
we say that condition (C) is satisfied if
there exists a constant $C$ such that
$$\int_0^1\omega_2(s/t)\,\frac{\psi(t)}{t}\,dt\le C\omega_1(s)$$
for every $s\in\Rp$.

% [{\bf ??} Probably not needed: 
% This is equivalent to the existence of a constant $C$ such that
% $$\int_s^{\infty}\frac{\omega_2(y)}{y}\,\psi(s/y)\,dy\le C\omega_1(s)$$
% for every $s\in\Rp$.]

\begin{theorem}
\label{th:char}
Let $\psi$ be a non-negative, measurable function on $[0,1]$ and let 
$\omega_1$ and $\omega_2$ be positive, continuous functions on $\Rp$.
Then $U_{\psi}$ defines a bounded operator from $\Loo$ to $\Lot$ if and only if
condition (C) is satisfied.
\end{theorem}

\noindent{\bf Proof}\quad 
Assume that condition (C) is satisfied and let $f\in\Loo$. Then
$$\int_0^{\infty}\int_0^1|f(s)|\,\frac{\psi(t)}{t}\,\omega_2(s/t)\,dt\,ds
\le C\int_0^{\infty}|f(s)|\omega_1(s)\,ds=C\|f\|_{\Loo}<\infty,$$
so it follows from Fubini's theorem that
$$\int_0^1\int_0^{\infty}|f(tx)|\psi(t)\omega_2(x)\,dx\,dt
=\int_0^1\int_0^{\infty}|f(s)|\,\frac{\psi(t)}{t}\,\omega_2(s/t)\,ds\,dt
\le C\|f\|_{\Loo}<\infty.$$
Another application of Fubini's theorem thus shows that 
$(U_{\psi}f)(x)$ is defined for almost all $x\in\Rp$ with
\begin{align*}
\|U_{\psi}f\|_{\Lot}
&=\int_0^{\infty}|(U_{\psi}f)(x)|\omega_2(x)\,dx  
\le\int_0^{\infty}\int_0^1|f(tx)|\psi(t)\omega_2(x)\,dt\,dx\\ 
&= \int_0^1\int_0^{\infty}|f(tx)|\psi(t)\omega_2(x)\,dx\,dt 
\le C\|f\|_{\Loo}<\infty.
\end{align*}
Hence $U_{\psi}$ defines a bounded operator from $\Loo$ to $\Lot$.
\medskip

Conversely, assume that $U_{\psi}$ defines a bounded operator from $\Loo$ to 
$\Lot$.
Since $\Lot$ is a closed subspace of $\Mot$ which we identify with the dual 
space of $\Cot$, 
it follows from \cite[Theorem~VI.8.6]{Du-Sc:I} that 
there exists a map $\rho$ from $\Rp$ to $\Mot$ for which the map
$s\mapsto\langle g,\rho(s)\rangle$ is measurable and essentially bounded
on $\Rp$ for every $g\in\Cot$ with 
$\|U_{\psi}\|=\text{ess\,sup}_{s\in\Rp}\|\rho(s)\|_{\Mot}$
and such that
$$\langle g,U_{\psi}f\rangle
=\int_0^{\infty}\langle g,\rho(s)\rangle f(s)\omega_1(s)\,ds
=\int_0^{\infty}\int_{\Rp}g(x)\,d\rho(s)(x)\,f(s)\omega_1(s)\,ds$$
for every $g\in\Cot$ and $f\in\Loo$. 
On the other hand
\begin{align*}
\langle g,U_{\psi}f\rangle
&=\int_0^{\infty}g(x)(U_{\psi}f)(x)\,dx\\
&=\int_0^{\infty}\int_0^x\frac{g(x)}{x}\,f(s)\psi(s/x)\,ds\,dx\\
&=\int_0^{\infty}\frac{1}{\omega_1(s)}
\int_s^{\infty}\frac{g(x)}{x}\,\psi(s/x)\,dx\,f(s)\omega_1(s)\,ds
\end{align*}
for every $g\in\Cot$ and $f\in\Loo$, so it follows that
$$\int_{\Rp}g(x)\,d\rho(s)(x)
=\frac{1}{\omega_1(s)}\int_s^{\infty}\frac{g(x)}{x}\,\psi(s/x)\,dx$$
for almost all $s\in\Rp$ and every $g\in\Cot$
(considering both sides as elements of $L^{\infty}(\Rp)$). 
Considered as elements of $\Mot$ we thus have
$$d\rho(s)(x)=\frac{1}{\omega_1(s)}\,\frac{1}{x}\,\psi(s/x)\,1_{x\ge s}\,dx$$
for almost all $s,x\in\Rp$. 
Hence $\rho(s)\in\Lot$ with
\begin{align*}
\|\rho(s)\|_{\Lot}
&=\int_0^{\infty}\omega_2(x)\,d\rho(s)(x)\\
&=\frac{1}{\omega_1(s)}\int_0^{\infty}\frac{1}{x}\,\psi(s/x)\,1_{x\ge s}\,\omega_2(x)\,dx
\\
&=\frac{1}{\omega_1(s)}\int_s^{\infty}\frac{1}{x}\,\psi(s/x)\omega_2(x)\,dx\\
&=\frac{1}{\omega_1(s)}\int_0^1\frac{\psi(t)}{t}\,\omega_2(s/t)\,dt
\end{align*}
for almost all $s\in\Rp$. Therefore
$$\int_0^1\omega_2(s/t)\,\frac{\psi(t)}{t}\,dt
=\|\rho(s)\|_{\Lot}\omega_1(s)\le \|U_{\psi}\|\omega_1(s)$$
for almost all $s\in\Rp$. 
Since both sides of the inequality are continuous functions of $s$, the 
inequality holds for every $s\in\Rp$, so condition (C) holds.
{\nopagebreak\hfill\raggedleft$\Box$\bigskip}

Letting $s=0$ in condition (C) we see that Xiao's condition is necessary
in our situation.

\begin{corollary}
\label{co:char}
Let $\psi$ be a non-negative, measurable function on $[0,1]$ and let 
$\omega_1$ and $\omega_2$ be positive, continuous functions on $\Rp$.
If $U_{\psi}$ defines a bounded operator from $\Loo$ to $\Lot$, then
$$\int_0^1\frac{\psi(t)}{t}\,dt<\infty.$$
\end{corollary}

The following straightforward consequences can be deduced from 
Theorem~\ref{th:char}.

\begin{corollary}
\label{co:char2}
Let $\psi$ be a non-negative, measurable function on $[0,1]$
\begin{enumerate}[(a)]
\item
Let $\omega$ be a decreasing, positive, continuous function on $\Rp$,
and assume that\\ 
$\int_0^1\psi(t)/t\,dt<\infty$.
Then $U_{\psi}$ defines a bounded operator from $\Lom$ to $\Lom$.
\item 
Let $\omega_1$ and $\omega_2$ be positive, continuous functions on $\Rp$,
and assume that $\omega_2$ is increasing. If $U_{\psi}$ defines a bounded 
operator from $\Loo$ to $\Lot$, then there exists a constant $C$ such that
$\omega_2(s)\le C\omega_1(s)$ for every $s\in\Rp$.
\item 
Let $\omega$ be an increasing, positive, continuous function on $\Rp$,
and assume that there exists $a<1$ and $K>0$ such that $\psi(t)\ge K$ 
almost everywhere on $[a,1]$. 
If $U_{\psi}$ defines a bounded operator from $\Lom$ to $\Lom$, then there 
exist positive constants $C_1$ and $C_2$ such that
$$C_1\omega(s)\le\int_0^1\omega(s/t)\,\frac{\psi(t)}{t}\,dt
\le C_2\omega(s)$$ 
for every $s\in\Rp$.
\end{enumerate}
\end{corollary}

\noindent{\bf Proof}\quad 
(a):\quad 
We have
$$\int_0^1\omega(s/t)\,\frac{\psi(t)}{t}\,dt
\le\int_0^1\frac{\psi(t)}{t}\,dt\,\omega(s)$$
for every $s\in\Rp$, so condition (C) is satisfied with 
$\omega_1=\omega_2=\omega$ and the result follows.

(b):\quad 
We have
$$\int_0^1\omega_2(s/t)\,\frac{\psi(t)}{t}\,dt
\ge\int_0^1\frac{\psi(t)}{t}\,dt\,\omega_2(s)$$
for every $s\in\Rp$. Since condition (C) is satisfied, the result follows.

(c):\quad 
We have
$$\int_0^1\omega(s/t)\,\frac{\psi(t)}{t}\,dt
\ge K\int_a^1\omega(s/t)\,dt\ge K(1-a)\omega(s)$$
for every $s\in\Rp$.
The other inequality is just condition (C) with 
$\omega_1=\omega_2=\omega$.
{\nopagebreak\hfill\raggedleft$\Box$\bigskip}

We finish the section with some examples of functions $\psi,\omega_1$ and 
$\omega_2$ for which $U_{\psi}$ defines a bounded operator from $\Loo$ to $\Lot$.

\begin{example}
\label{ex:char}
\ 
\begin{enumerate}[(a)]
\item
For $\alpha>0$ and $\beta_1,\beta_2\in\mathbb R$, let $\psi(t)=t^{\alpha}$ for $t\in[0,1]$
and $\omega_i(x)=(1+x)^{\beta_i}$ for $x\in\Rp$ and $i=1,2$. 
Then $U_{\psi}$ defines a bounded operator from $\Loo$ to $\Lot$ if and only if
$\beta_2\le\beta_1$ and $\beta_2<\alpha$.
\item 
For $\alpha>0$, let $\psi(t)=t^{\alpha}$ for $t\in[0,1]$. Also,
let $\omega_1(x)=e^{-x}/(1+x)$ and $\omega_2(x)=e^{-x}$ for $x\in\Rp$.
Then $U_{\psi}$ defines a bounded operator from $\Loo$ to $\Lot$.
Moreover, it is not possible to replace $\omega_1(x)$ by a function tending faster to zero at 
infinity.
\item 
Let $\psi(t)=e^{-1/t^2}$ for $t\in[0,1]$. Also, let $\omega_1(x)=e^{x^2/4}/x$ and 
$\omega_2(x)=e^{x}$ for $x\in\Rp$.
Then $U_{\psi}$ defines a bounded operator from $\Loo$ to $\Lot$.
Moreover, it is not possible to replace $\omega_1(x)$ by a function tending slower to infinity 
at infinity.
\end{enumerate}
\end{example}

\noindent{\bf Proof}\quad 
(a):\quad 
For $s\ge1$ and $t\in[0,1]$ we have $s/t<1+s/t\le2s/t$, so
\begin{align*}
\int_0^1\omega_2(s/t)\,\frac{\psi(t)}{t}\,dt
&= \int_0^1\left(1+\frac{s}{t}\right)^{\beta_2}t^{\alpha-1}\,dt\\
&\simeq s^{\beta_2}\int_0^1t^{\alpha-\beta_2-1}\,dt\\
&\simeq s^{\beta_2}
\end{align*}
for $s\ge1$ if $\beta_2<\alpha$, whereas the integrals diverge if $\beta_2\ge\alpha$. 
Moreover, the expression 
$$\int_0^1\omega_2(s/t)\,\frac{\psi(t)}{t}\,dt
=\int_0^1\left(1+\frac{s}{t}\right)^{\beta_2}t^{\alpha-1}\,dt$$
defines a positive, continuous function of $s$ on $\Rp$, so it follows that  
condition (C) is satisfied if and only if $\beta_2\le\beta_1$ and $\beta_2<\alpha$.

(b):\quad 
For $s\ge1$ we have 
$$\int_0^1\omega_2(s/t)\,\frac{\psi(t)}{t}\,dt
= \int_s^{\infty}\frac{\omega_2(x)}{x}\,\psi(s/x)\,dx
= \int_s^{\infty}\frac{e^{-x}}{x}\,\frac{s^{\alpha}}{x^{\alpha}}\,dx
\le \int_s^{\infty}\frac{e^{-x}}{x}\,dx
\le \frac{e^{-s}}{s}\,.$$
Moreover, 
$$\int_0^1\omega_2(s/t)\,\frac{\psi(t)}{t}\,dt
\le\int_0^1\frac{\psi(t)}{t}\,dt<\infty$$ 
for all $s\in\Rp$, so condition (C) is satisfied and $U_{\psi}$ thus defines a 
bounded operator from $\Loo$ to $\Lot$.
On the other hand, since
$$\int_0^1\omega_2(s/t)\,\frac{\psi(t)}{t}\,dt
% = \int_s^{\infty}\frac{e^{-x}}{x}\,\frac{s^{\alpha}}{x^{\alpha}}\,dx
\ge \int_s^{2s}\frac{e^{-x}}{x}\,\frac{s^{\alpha}}{x^{\alpha}}\,dx
\ge \frac{1}{2^{\alpha+1}s}\int_s^{2s}e^{-x}\,dx
\ge \frac{1}{2^{\alpha+2}}\,\frac{e^{-s}}{s}$$
for $s\ge1$, it is not possible to replace $\omega_1(x)$ by a function tending 
faster to zero at infinity.

(c):\quad 
For $s\in\Rp$ we have
$$\int_0^1\omega_2(s/t)\,\frac{\psi(t)}{t}\,dt
= \int_s^{\infty}\frac{\omega_2(x)}{x}\,\psi(s/x)\,dx
= \int_s^{\infty}\frac{e^{x-x^2/s^2}}{x}\,dx
= \int_1^{\infty}\frac{e^{sy-y^2}}{y}\,dy.$$
Moreover, for $s\ge4$ 
$$\int_{s/4}^{\infty}\frac{e^{sy-y^2}}{y}\,dy
\le\frac{4}{s}\int_{s/4}^{\infty}e^{-(y-s/2)^2+s^2/4}\,dy
=4\int_{-s/4}^{\infty}e^{-u^2}\,du\,\frac{e^{s^2/4}}{s}$$
and
$$\int_1^{s/4}\frac{e^{sy-y^2}}{y}\,dy
\le\int_1^{s/4}e^{sy}\,dy
\le\frac{e^{s^2/4}}{s},$$
so condition (C) is satisfied and $U_{\psi}$ thus defines a bounded operator from $\Loo$ to 
$\Lot$. 
On the other hand, the estimate 
$$\int_0^1\omega_2(s/t)\,\frac{\psi(t)}{t}\,dt
=\int_1^{\infty}\frac{e^{sy-y^2}}{y}\,dy
\ge\frac{1}{s}\int_{s/2}^{s/2+1}e^{-(y-s/2)^2+s^2/4}\,dy
=\int_{0}^{1}e^{-u^2}\,du\,\frac{e^{s^2/4}}{s}$$
for $s\ge2$ shows that it is not possible to replace $\omega_1(x)$ by a 
function tending slower to infinity at infinity.
{\nopagebreak\hfill\raggedleft$\Box$\bigskip}

In Example~\ref{ex:char}(b) we have $\omega_2(x)/\omega_1(x)\to\infty$ as $x\to\infty$,
which should be compared to the conclusion in Corollary~\ref{co:char2}(b). 
Conversely, Example~\ref{ex:char}(c) shows an example where we need 
$\omega_2(x)/\omega_1(x)\to0$ rapidly as $x\to\infty$ in order for $U_{\psi}$ to be defined.

\section{Extensions to weighted spaces of measures}
\label{sec:ext}

Identifying the dual space of $\Lom$ with $\Lim$ as in the introduction,
we have the following result about the adjoint of $U_{\psi}$.

\begin{proposition}
\label{pr:dual}
Let $\psi$ be a non-negative, measurable function on $[0,1]$ and let 
$\omega_1$ and $\omega_2$ be positive, continuous functions on $\Rp$.
Assume that condition (C) is satisfied so that 
$U_{\psi}:\Loo\to\Lot$ is a bounded operator, and consider the adjoint operator
$U_{\psi}^*:\Lit\to\Lio$. 
\begin{enumerate}[(a)]
\item
For $h\in\Lit$ we have
$$(U_{\psi}^*h)(x)=\int_0^1 h(x/t)\,\frac{\psi(t)}{t}\,dt$$
% \quad\text{[?? include here]}=\int_x^\infty h(s)\,\frac{1}{x}\,\psi(x/s)\,ds$$
for almost all $x\in\Rp$.
\item 
$U_{\psi}^*$ maps $\Cot$ into $\Coo$.
\end{enumerate}
\end{proposition}

\noindent{\bf Proof}\quad 
(a):\quad 
Let $h\in\Lit$. Since $|h(x/t)|\le\|h\|_{\Lit}\omega_2(x/t)$ for almost all 
$x,t\in\Rp$, it follows from condition (C) that 
$\int_0^1 h(x/t)\psi(t)/t\,dt$ is defined and satisfies
$$\left|\int_0^1 h(x/t)\,\frac{\psi(t)}{t}\,dt\right|
\le\|h\|_{\Lit}\int_0^1 \omega_2(x/t)\,\frac{\psi(t)}{t}\,dt
\le C\|h\|_{\Lit}\omega_1(x)$$
for almost all $x\in\Rp$. 
Hence the function $x\mapsto\int_0^1 h(x/t)\psi(t)/t\,dt$ 
belongs to $\Lio$. Also, for $f\in\Loo$ we have
\begin{align*}
\langle f,U_{\psi}^*h\rangle=\langle U_{\psi}f,h\rangle
&=\int_0^{\infty}(U_{\psi}f)(s)h(s)\,ds\\
&=\int_0^{\infty}\int_0^s\frac{1}{s}\,f(x)\psi(x/s)h(s)\,dx\,ds\\
&=\int_0^{\infty}\int_x^{\infty}\frac{h(s)}{s}\,\psi(x/s)\,ds\,f(x)\,dx
\end{align*}
from which it follows that
$$(U_{\psi}^*h)(x)=\int_x^{\infty}\frac{h(s)}{s}\,\psi(x/s)\,ds
=\int_0^1 h(x/t)\,\frac{\psi(t)}{t}\,dt$$
for almost all $x\in\Rp$.

(b):\quad
It suffices to show that $U_{\psi}^*$ maps $C_c(\Rp)$ 
(the continuous functions on $\Rp$ with compact support) into $\Coo$.
Let $g\in C_c(\Rp)$, let $x_0\in\Rp$ and let $(x_n)$ be a sequence in $\Rp$ 
with $x_n\to x_0$ as $n\to\infty$. Then 
$$(U_{\psi}^*g)(x_n)-(U_{\psi}^*g)(x_0)
=\int_0^1 (g(x_n/t)-g(x_0/t))\,\frac{\psi(t)}{t}\,dt$$
for $n\in\mathbb N$.
Since $g$ is bounded on $\Rp$ and since $\int_0^1\psi(t)/t\,dt<\infty$ 
by Corollary~\ref{co:char}, it follows from Lebesgue's dominated convergence 
theorem that $(U_{\psi}^*g)(x_n)\to(U_{\psi}^*g)(x_0)$ as $n\to\infty$. Hence 
$U_{\psi}^*g$ is continuous on $\Rp$.
Finally, from the expression
$$(U_{\psi}^*g)(x)=\int_x^{\infty}\frac{g(s)}{s}\,\psi(x/s)\,ds$$
it follows that $\text{supp}\,U_{\psi}^*g\subseteq\text{supp}\,g$, so we
conclude that $U_{\psi}^*g\in C_c(\Rp)\subseteq\Coo$.
{\nopagebreak\hfill\raggedleft$\Box$\bigskip}

Let $V_{\psi}$ be the restriction of $U_{\psi}^*$ to $\Cot$ considered as a map
into $\Coo$. We then immediately have the following result.

\begin{corollary}
\label{co:ext}
Let $\psi$ be a non-negative, measurable function on $[0,1]$ and let 
$\omega_1$ and $\omega_2$ be positive, continuous functions on $\Rp$.
Assume that condition (C) is satisfied so that 
$U_{\psi}:\Loo\to\Lot$ is a bounded operator. 
The bounded operator $\overline{U}_{\psi}=V_{\psi}^*$ from $\Moo$ to $\Mot$ is an 
extension of $U_{\psi}$.
\end{corollary}

Let $\psi$ be a non-negative, continuous function on $[0,1]$ with $\psi(0)=0$. 
For $\mu\in\Moo$ and $x>0$ let 
$$(W_{\psi}\mu)(x)=\frac{1}{x}\int_{(0,x)}\psi(s/x)\,d\mu(s).$$

\begin{proposition}
\label{pr:ext1}
Let $\psi$ be a non-negative, continuous function on $[0,1]$ and let 
$\omega_1$ and $\omega_2$ be positive, continuous functions on $\Rp$.
Assume that condition (C) is satisfied so that 
$U_{\psi}:\Loo\to\Lot$ is a bounded operator. 
Then $W_{\psi}\mu\in\Lot$ and
$$\overline{U}_{\psi}\mu
=W_{\psi}\mu+\int_0^1\frac{\psi(t)}{t}\,dt\cdot\mu(\{0\})\delta_0$$
for $\mu\in\Moo$. 
In particular $\text{\em{ran}}\,\overline{U}_{\psi}\subseteq\Lot\oplus\mathbb C\delta_0$
and $\overline{U}_{\psi}$ maps $M((0,\infty),\omega_1)$ into $\Lot$.
\end{proposition}

\noindent{\bf Proof}\quad 
By Corollary~\ref{co:char} we have $\int_0^1\psi(t)/t\,dt<\infty$, so it 
follows that $\psi(0)=0$. 
Let $\mu\in\Moo$ with $\mu(\{0\})=0$. By condition (C) we have
\begin{align*}
\int_{(0,\infty)}\int_s^{\infty}\frac{1}{x}\,\psi(s/x)\omega_2(x)
\,dx\,d|\mu|(s)
& =\int_{(0,\infty)}\int_0^1\omega_2(s/t)\,\frac{\psi(t)}{t}\,dt\,d|\mu|(s)\\
& \le C\int_{(0,\infty)}\omega_1(s)\,d|\mu|(s)=C\|\mu\|_{\Moo}<\infty,
\end{align*}
so it follows from Fubini's theorem that
$$\int_0^{\infty}\frac{1}{x}\int_{(0,x)}\,\psi(s/x)\,d|\mu|(s)\,\omega_2(x)\,dx
<\infty.$$
Hence $W_{\psi}\mu\in\Lot$. Moreover, for $g\in\Cot$ we have
\begin{align*}
\langle g,\overline{U}_{\psi}\mu\rangle=\langle V_{\psi}g,\mu\rangle
&= \int_{(0,\infty)}\int_0^1g(s/t)\,\frac{\psi(t)}{t}\,dt\,d\mu(s)\\
&= \int_{(0,\infty)}\int_s^{\infty}\frac{g(x)}{x}\,\psi(s/x)\,dx\,d\mu(s)\\
&= \int_0^{\infty}\frac{1}{x}\int_{(0,x)}\,\psi(s/x)\,d\mu(s)\,g(x)\,dx\\
&= \int_0^{\infty}(W_{\psi}\mu)(x)g(x)\,dx=\langle g,W_{\psi}\mu\rangle,
\end{align*}
so we conclude that $\overline{U}_{\psi}\mu=W_{\psi}\mu$. 
Finally, for $g\in\Cot$ we have
$$\langle g,\overline{U}_{\psi}\delta_0\rangle=\langle V_{\psi}g,\delta_0\rangle
=(V_{\psi}g)(0)=g(0)\int_0^1\frac{\psi(t)}{t}\,dt
=\langle g,\int_0^1\frac{\psi(t)}{t}\,dt\cdot\delta_0\rangle.$$
Since $W_{\psi}\delta_0=0$ this finishes the proof.
{\nopagebreak\hfill\raggedleft$\Box$\bigskip}

The conclusion about the range of $\overline{U}_{\psi}$ can be generalized
to the case, where $\psi$ is not assumed to be continuous.

\begin{proposition}
\label{pr:ext2}
Let $\psi$ be a non-negative, measurable function on $[0,1]$ and let 
$\omega_1$ and $\omega_2$ be positive, continuous functions on $\Rp$.
Assume that condition (C) is satisfied so that 
$U_{\psi}:\Loo\to\Lot$ is a bounded operator. 
Then $\text{\em{ran}}\,\overline{U}_{\psi}\subseteq\Lot\oplus\mathbb C\delta_0$.
\end{proposition}

\noindent{\bf Proof}\quad 
Choose a sequence of non-negative, continuous functions 
$(\psi_n)$ on $[0,1]$ with $\psi_n\le\psi$ and
$$\int_0^1\frac{\psi(t)-\psi_n(t)}{t}\,dt\to0\qquad\text{as }n\to\infty.$$
For $\mu\in\Moo$ and $g\in\Cot$ we have
\begin{align*}
|\langle g,(\overline{U}_{\psi}-\overline{U}_{\psi_n})\mu\rangle|
&=|\langle (V_{\psi}-V_{\psi_n})g,\mu\rangle|\\
&= \left|\int_{\Rp}\int_0^1g(x/t)\,\frac{\psi(t)-\psi_n(t)}{t}
\,dt\,d\mu(x)\right|\\
&\le\|g\|_{\Cot}\int_{\Rp}\int_0^1\omega_2(x/t)\,\frac{\psi(t)-\psi_n(t)}{t}
\,dt\,d|\mu|(x).
\end{align*}
Let 
$$p_n(x)=\int_0^1\omega_2(x/t)\,\frac{\psi(t)-\psi_n(t)}{t}\,dt$$
for $x\in\Rp$ and $n\in\mathbb N$. By condition (C) there exists a constant
$C$ such that $p_n(x)\le C\omega_1(x)$ for every $x\in\Rp$ and $n\in\mathbb N$.
Moreover, for every $x\in\Rp$ we have $p_n(x)\to0$ as $n\to\infty$ by 
Lebesgue's dominated convergence theorem. Hence 
$$\|(\overline{U}_{\psi}-\overline{U}_{\psi_n})\mu\|_{\Mot}=
\sup_{\|g\|_{\Cot}\le1}|\langle g,(\overline{U}_{\psi}-\overline{U}_{\psi_n})\mu\rangle|
\le\int_{\Rp}p_n(x)\,d|\mu|(x)\to0$$
as $n\to\infty$
again by Lebesgue's dominated convergence theorem. Consequently,
$\overline{U}_{\psi_n}\to\overline{U}_{\psi}$ strongly as $n\to\infty$.
Since $\text{ran}\,\overline{U}_{\psi_n}\subseteq\Lot\oplus\mathbb C\delta_0$ for
$n\in\mathbb N$ by Proposition~\ref{pr:ext1}, the same thus holds for
$\text{ran}\,\overline{U}_{\psi}$. 
{\nopagebreak\hfill\raggedleft$\Box$\bigskip}

\begin{corollary}
\label{co:ext2}
Let $\psi$ be a non-negative, measurable function on $[0,1]$ and let 
$\omega_1$ and $\omega_2$ be positive, continuous functions on $\Rp$.
Assume that condition (C) is satisfied so that 
$U_{\psi}:\Loo\to\Lot$ is a bounded operator. 
For $s>0$ we then have $(\overline{U}_{\psi}\delta_s)(x)=\psi(s/x)/x$ 
for almost all $x\ge s$ and $(\overline{U}_{\psi}\delta_s)(x)=0$ for almost all $x<s$.
\end{corollary}

\noindent{\bf Proof}\quad 
For $\psi$ continuous, this follows from Proposition~\ref{pr:ext1}. 
For general $\psi$ it follows from the approach in the proof of 
Proposition~\ref{pr:ext2} using $\overline{U}_{\psi_n}\to\overline{U}_{\psi}$ strongly as
$n\to\infty$.
{\nopagebreak\hfill\raggedleft$\Box$\bigskip}

It follows from Corollary~\ref{co:ext2} that
$$\|\overline{U}_{\psi}\delta_s\|_{\Mot}
=\int_s^{\infty}\frac{\omega_2(x)}{x}\,\psi(s/x)\,dx
=\int_0^1\omega_2(s/t)\,\frac{\psi(t)}{t}\,dt,$$
whereas $\|\delta_s\|_{\Moo}=\omega_1(s)$. 
Since $\overline{U}_{\psi}$ is bounded we thus recover condition (C).
If we without using Theorem~\ref{th:char} could show that if 
$U_{\psi}:\Loo\to\Lot$ is a bounded operator, then is has a bounded extension
$\overline{U}_{\psi}:\Moo\to\Mot$ for which Corollary~\ref{co:ext2} holds, then
we would in this way obtain an alternative proof of condition (C).

\section{Weakly compact operators}
\label{sec:wk}

We finish the paper by showing that there are no non-zero, weakly compact generalized 
Hardy-Ces\`aro operators between $\Loo$ and $\Lot$.

\begin{proposition}
\label{co:wc}
Let $\psi$ be a non-negative, measurable function on $[0,1]$ and let 
$\omega_1$ and $\omega_2$ be positive, continuous functions on $\Rp$.
Assume that condition (C) is satisfied so that 
$U_{\psi}:\Loo\to\Lot$ is a bounded operator. 
If $\psi\neq0$, then $U_{\psi}$ is not weakly compact.
\end{proposition}

\noindent{\bf Proof}\quad 
For $f\in\Loo$ and $x\in\Rp$ we have
$$(U_{\psi}f)(x)=\frac{1}{x}\int_0^x f(s)\psi(s/x)\,ds
=\int_0^{\infty}f(s)\rho(s)(x)\omega_1(s)\,ds,$$
where (with a slight change of notation compared to the proof of 
Theorem~\ref{th:char}) 
$$\rho(s)(x)=\frac{1}{\omega_1(s)}\,\frac{1}{x}\,\psi(s/x)\,1_{x\ge s}$$
for $x,s\in\Rp$. In the proof of Theorem~\ref{th:char} we saw that
$\rho(s)\in\Lot$ with $\|\rho(s)\|_{\Lot}\le C$ for a constant $C$ for almost all
$s\in\Rp$.
It thus follows from \cite[Theorem~VI.8.10]{Du-Sc:I} that 
$U_{\psi}$ is weakly compact if and only if $\{\rho(s):s\in\Rp\}$
is contained in a weakly compact set of $\Lot$ 
(except possibly for $s$ belonging to a null-set). 
Consider $\rho(s)$ as an element of 
$\Cot^*$ for $s\in\Rp$ and let $g\in\Cot$. Then
\begin{align*}
\langle g,\rho(s)\rangle
&=\int_0^{\infty}g(x)\rho(s)(x)\,dx\\
&=\frac{1}{\omega_1(s)}\int_s^{\infty}\frac{g(x)}{x}\,\psi(s/x)\,dx\\
&=\frac{1}{\omega_1(s)}\int_0^1g(s/t)\,\frac{\psi(t)}{t}\,dt.
\end{align*}
Since $g(s/t)\to g(0)$ as $s\to0_+$ for all $t>0$, it follows from
Lebesgue's dominated convergence theorem that
$$\langle g,\rho(s)\rangle\to
\frac{1}{\omega_1(0)}\,g(0)\int_0^1\frac{\psi(t)}{t}\,dt$$
as $s\to0_+$. We therefore conclude that
$$\rho(s)\to\frac{1}{\omega_1(0)}\int_0^1\frac{\psi(t)}{t}\,dt\cdot\delta_0$$
weak-star in $\Mot$ as $s\to0_+$. Since $\delta_0\notin\Lot$, it follows that 
$\{\rho(s):s\in\Rp\}$
is not contained in a weakly compact set of $\Lot$ (even excepting null sets),
and the result follows.
{\nopagebreak\hfill\raggedleft$\Box$\bigskip}

\bigskip
\bigskip

% \bibliography{myref}

% \bibliographystyle{plain}

\bigskip

\noindent
Thomas Vils Pedersen\\
Department of Mathematics\\
University of Copenhagen\\
Universitetsparken 5\\
DK-2100 Copenhagen \O\\
Denmark\\
vils@math.ku.dk

\end{document}